\documentclass{article}

\usepackage{amsmath,amsthm,amsfonts,amssymb,amsopn}

\title{A Hecke Correspondence Theorem for Automorphic Integrals with Symmetric Rational Period Functions on the Hecke Groups}

\author{Wendell Ressler \\
        Franklin \& Marshall College \\
        Lancaster, PA 17604}

\date{August 7, 2008}

   \theoremstyle{plain} 					
   \newtheorem{theorem}{Theorem}
   
   \newtheorem{lemma}{Lemma}

   \theoremstyle{definition}

   \theoremstyle{remark}
   
   \newtheorem*{remarks}{Remarks}

   \newcommand{\Z}{\mathbb{Z}}		
   \newcommand{\R}{\mathbb{R}}		
   \newcommand{\C}{\mathbb{C}}		
   
   \newcommand{\re}{\text{Re}}		
   
   \newcommand{\uhp}{\mathcal{H}}	
   \newcommand{\Pset}{\mathcal{P}}	
   
   \DeclareMathOperator*{\res}{Res}	

\begin{document}

\maketitle

\section{Introduction}

In the 1930s Hecke \cite{Hec36, Hec38} formalized a general correspondence between automorphic forms and Dirichlet series.
Hecke's work generalized Riemann's use of the transformation law for the elliptic
\( \theta \)-function to derive the functional equation for the zeta function
\( \zeta(s) \) in
\cite{Rie1859}.

Hecke studied automorphic forms with respect to an infinite class of 
discrete groups that act 
on the upper half plane as linear 
fractional transformations.
These groups 
have become known as the Hecke groups,
and include the modular group
\( \Gamma(1) = \mathrm{PSL}(2,\Z) \).

In \cite{Eic57} Eichler introduced generalized abelian integrals, which he obtained by integrating modular forms of positive weight.
An Eichler integral satisfies a modular relation with a polynomial period function.
In \cite{Kno78} and \cite{Kno81} Knopp generalizes Eichler integrals
and develops the theory of automorphic integrals with rational period functions.
Knopp shows that an entire modular integral corresponds to a Dirichlet
series that satisfies Hecke's functional equation,
provided the rational period function has poles only at 
\( 0 \) or \( \infty \).
Knopp also proves a converse theorem, from which it follows that if a rational period function has any other poles the corresponding Dirichlet
series cannot satisfy the same functional equation.

In \cite{HaK92} Hawkins and Knopp prove a Hecke correspondence theorem for modular integrals with rational period functions on 
\( \Gamma_{\theta} \), the theta subgroup of 
\( \Gamma(1) \).  
In this correspondence the Dirichlet series functional equation contains a remainder term that corresponds to the nonzero poles of the rational period function.
Hawkins and Knopp observe that their theorem implies that an automorphic integral with a rational period function on one of the Hecke groups must correspond to a Dirichlet series that satisfies a functional equation similar to the one they found.
Hawkins and Knopp also point out that since the Hecke groups have two group relations,
while \( \Gamma_{\theta} \) has a single group relation,
rational period functions on Hecke groups have more structure than rational period functions on 
\( \Gamma_{\theta} \).
Thus for Hecke groups the corresponding remainder terms must have more structure than the ones discovered by Hawkins and Knopp,
a fact that a full correspondence theorem in this setting must reveal.

In \cite{CR96} this author proves a Hecke correspondence theorem for modular integrals with rational period functions on
\( \Gamma(1) \), which is one of the Hecke groups.  
We show that the remainder term for the 
Dirichlet series functional equation satisfies a second relation
that corresponds to the second group relation in 
\( \Gamma(1) \).

In this paper we extend the correspondence to a class of
automorphic integrals with rational period functions 
on \emph{any} of the Hecke groups.  
We restrict our attention to automorphic integrals 
of weight that is twice an odd integer
and to rational period functions that satisfy
a certain symmetry property we call 
``Hecke-symmetry.''
We show that the remainder term in the 
Dirichlet series functional equation
satisfies a second relation that 
generalizes the second relation in \cite{CR96}.

\section{Hecke groups and fixed points}

Let \( \lambda \) be a fixed positive real number and
put \(S = S_{\lambda} 
   = \left( \begin{smallmatrix} 1&\lambda\\ 0&1 \end{smallmatrix} \right) \),
\(T = \left( \begin{smallmatrix} 0&-1\\ 1&0 \end{smallmatrix} \right)\),
and \(I = \left( \begin{smallmatrix} 1&0\\ 0&1 \end{smallmatrix} \right)\).
Define the group
\(G(\lambda)=\langle S,T \rangle / \{\pm I\} \subseteq \mathrm{PSL}(2,\R) \).
Elements of this group
act on the Riemann sphere as linear fractional transformations,
that is,  
\( Mz = \frac{az+b}{cz+d} \) for
\( M = \left( \begin{smallmatrix} a&b\\ c&d \end{smallmatrix} \right) 
\in G(\lambda)\)
and \( z \in \C \cup \left\{ \infty \right\} \).
This action preserves the real line and the upper half-plane
\( \mathcal{H} \).

Erich Hecke \cite{Hec36,Hec38} showed that the 
values of \(\lambda\) between \( 0 \) and \( 2 \)
for which \(G(\lambda)\) is discrete 
are 
\begin{equation*}
    \lambda = \lambda_{p} = 2\cos(\pi/p),
\end{equation*}
for \(p=3, 4, 5, \dots\).
These discrete groups are the 
\emph{Hecke groups}, 
which we denote by
\(G_{p}=G(\lambda_{p})\) for \(p\geq3\).
Each of the Hecke groups has two group relations,
which may be written
\( T^{2}=(ST)^{p}=I \).
(Note that we are identifying
\( I \) and \( -I \),
since these are projective groups.)
The first of these groups is the modular group
\(G_{3}=G(1)=\Gamma(1)\).
For the rest of the paper we fix the integer
\( p \geq 3 \) and the real number \( \lambda = \lambda_{p} \).

An element \(M = \bigl( \begin{smallmatrix} a&b\\ c&d \end{smallmatrix} \bigr)
    \in G_{p}\) is 
\emph{hyperbolic} if \( \vert a+d \vert >2 \),
\emph{parabolic} if \( \vert a+d \vert =2 \),
and \emph{elliptic} if \( \vert a+d \vert <2 \). 
We designate fixed points accordingly.
The element
\( M = \bigl( \begin{smallmatrix} a&b\\ c&d \end{smallmatrix} \bigr) \)
fixes
\begin{eqnarray}
    z & = & \frac{a-d \pm \sqrt{(d-a)^{2}+4bc}}{2c} \nonumber \\
      & = & \frac{a-d \pm \sqrt{(a+d)^{2}-4}}{2c},
    \label{FixedPointFormula}
\end{eqnarray}
so hyperbolic elements of 
\( G_{p} \) have two distinct real fixed points.
Since \( G_{p} \) is discrete,
the stabilizer of any complex number
is a cyclic subgroup of \( G_{p} \)
\cite[page 15]{Leh66}.
If \( \alpha \) is a hyperbolic fixed point of 
\( G_{p} \) we call the 
other point fixed by its stabilizer the 
\emph{Hecke conjugate of}
\( \alpha \),
and we denote it by 
\( \alpha^{\prime} \).
A straightforward calculation shows that if 
\( \alpha \) is hyperbolic and
\( M \in G_{p} \),
then \( \left(M\alpha\right)^{\prime} = M\alpha^{\prime} \).
If \( R \) is a set of hyperbolic fixed points of 
\( G_{p} \) we write
\( R^{\prime} = \{x^{\prime} \mid x \in R\} \).
We say that \( R \) has 
\emph{Hecke symmetry}
if \( R = R^{\prime} \).

\section{Automorphic integrals}

Suppose \( F \) is a function holomorphic in 
the upper half-plane
\( \cal H \) with the Fourier expansion
\begin{equation}
	\label{Fourierexp}
	F(z) = \sum _{n=0}^{\infty} a_{n} e^{2 \pi inz/\lambda},
\end{equation}
for \( z \in \uhp \).
If for every $z \in \cal H$, 
$F$ satisfies the automorphic relation 
\begin{equation}
	\label{mirelation}
	z^{-2k} F\left( \frac{-1}{z} \right) = F(z) + q(z),		
\end{equation}
where \( q(z) \) is a rational function
and \( 2k \in 2\Z^{+} \),
we say that $F$ is an 
{\em entire automorphic integral of weight \( 2k \) on \( G_{p} \)
with rational period function (RPF) \( q \)}.
If \( q \equiv 0 \) then
\( F \) is an entire automorphic \emph{form}
of weight \( 2k \) on \( G_{p} \).
For 
\(M = \bigl( \begin{smallmatrix} *&*\\ c&d \end{smallmatrix} \bigr) 
\in G_{p}\)
and \(F(z)\) a complex function,
we define the
\emph{weight \( 2k \) slash operator}
\(F \mid_{2k} M = F \mid M\) by
\begin{equation*}
	\left( F\mid M \right) (z) = (cz+d)^{-2k} F \left( Mz \right).
\end{equation*}
With this notation we may rewrite the automorphic relation
\eqref{mirelation} as
\begin{equation}
	\label{}
	F \mid T = F + q.
	\nonumber
\end{equation}

A calculation shows that 
\( F \mid M_{1}M_{2} = \left( F \mid M_{1} \right) \mid M_{2} \)
for \( M_{1}, M_{2} \in G_{p} \).
We use this and the group relation 
\( T^{2} = I \) to calculate that a RPF 
\( q \) satisfies the relation
\begin{equation}
	\label{rpfrelation1}
	q \mid T + q = 0.
\end{equation}
In a similar way the relation 
\( (ST)^{p} = I \) implies that \( q \) satisfies a second relation
\begin{equation}
	\label{rpfrelation2}
	q + q \mid ST + \cdots + q \mid (ST)^{p-1} = 0.
\end{equation}
Knopp \cite[Section II]{Kno74} showed that 
(\ref{rpfrelation1}) and (\ref{rpfrelation2}) 
characterize the set of RPFs for a given 
group and weight.

\section{Binary quadratic forms}

Hawkins \cite{Haw} pointed out a deep connection between RPFs
on the modular group and classical
binary quadratic forms.
Schmidt \cite{Sch93} and Schmidt and Sheingorn \cite{SS95} 
observed that similar connections exist between
RPFs on the Hecke groups
and binary quadratic forms with coefficients in
\( \Z \left[ \lambda_{p} \right] \).
We give properties of these 
binary quadratic forms in \cite{Res08}
and we exploit the connections to RPFs 
on the Hecke groups
in \cite{CR01}.
In this section we describe the results we need
to describe the structure of RPFs and to
develop our correspondence.

We consider binary quadratic forms
\[ Q(x,y) = Ax^{2}+Bxy+Cy^{2}, \]
with coefficients in  
\( \Z[\lambda] \).
We denote such a form by \( Q=[A,B,C] \) and refer to it as a
\( \lambda_{p} \)-BQF or a \( \lambda \)-BQF.
We restrict our attention to indefinite forms,
which have positive discriminant 
\( D=B^{2}-4AC \).

Elements of a Hecke group act on \( \lambda \)-BQFs by
\( \left(Q \circ M \right)(x,y) = Q(a x + b y, c x + d y) \)
for 
\( Q \) a \( \lambda \)-BQF
and
\( M = \bigl( \begin{smallmatrix} a &b \\ c & d 
	      \end{smallmatrix} \bigr)\ \in G_{p} \).
Since
\( Q \circ M = Q \circ (-M) \), 
this definition does not depend on our 
choice of coset representative.
A calculation shows that 
\( B^{\prime 2}-4A^{\prime}C^{\prime} = B^{2}-4AC \),
so the action of a Hecke group preserves the discriminant.
We say that \( Q \) and \( \tilde{Q} \) are
\emph{\( G_{p} \)-equivalent}, and write \( Q \sim \tilde{Q} \),
if there exists a \( V \in G_{p} \) such that 
\( \tilde{Q} = Q \circ V \).
\( G_{p} \)-equivalence is an equivalence relation, so 
\( G_{p} \) partitions the set of \( \lambda \)-BQFs into equivalence 
classes of forms.

In \cite{Res08}
we describe a one-to-one correspondence between hyperbolic fixed points
of \( G_{p} \) and certain 
\( \lambda_{p} \)-BQFs.
In order to do so we use a variant of Rosen's 
\( \lambda \)-continued fractions 
\cite{Ros54, Sch93}
to first map every hyperbolic point 
\( \alpha \) to the unique primitive hyperbolic element
\( M_{\alpha} \in G_{p} \) with positive trace
that has \( \alpha \) as an attracting fixed point.
This mapping associates Hecke conjugates with inverse elements in the Hecke group, 
that is,
\( M_{\alpha^{\prime}} = M_{\alpha}^{-1} \).
We then map every primitive hyperbolic element
\( M = \left( \begin{smallmatrix} a&b\\ c&d \end{smallmatrix} \right) 
\in G_{p} \) 
with positive trace to the 
\( \lambda \)-BQF 
\( Q = [c,d-a,-b] \).
The roots of
\( Q(z,1) = cz^{2} + (d-a)z - b \)
are the
fixed points
of \( M \) given by
\eqref{FixedPointFormula}.
Since every domain element for this map
is hyperbolic 
we call the images \emph{hyperbolic} \( \lambda \)-BQFs.
We also note that every hyperbolic
\( \lambda \)-BQF
is indefinite.
The composition of the two maps above associates every hyperbolic fixed point
\( \alpha \) with a unique hyperbolic
\( \lambda \)-BQF that we denote
\( Q_{\alpha} \).
These mappings are all injective, so
the inverse of the composition exists;
this inverse associates every 
hyperbolic
quadratic form
\( Q = [A,B,C] \) with the hyperbolic number
\( \alpha_{Q} = \frac{-B+\sqrt{D}}{2A} \),
where 
\( D \) is the discriminant of \( Q \).

Suppose that \( V \in G_{p} \)
and 
\( \alpha \) and \( \beta \) are hyperbolic numbers.
Then 
\( \beta = V^{-1}\alpha \)
if and only if 
\( Q_{\beta} = Q_{\alpha} \circ V \)
for the associated forms, and
if and only if
\( M_{\beta} = V^{-1}M_{\alpha}V \)
for the associated matrices.
Now trace is preserved by conjugation, as is 
the property of a matrix being primitive.
Thus every \( G_{p} \)-equivalence class of
\( \lambda \)-BQFs contains either only hyperbolic forms
or no hyperbolic forms,
so we designate \( \lambda \)-BQF
equivalence classes themselves as hyperbolic or not hyperbolic.

If a hyperbolic quadratic form \( Q = [A,B,C] \) satisfies
\( A > 0 > C \)
we say that \( Q \) is
\( G_{p} \)-\emph{simple}
(or \emph{simple}, if the context is clear).
If \( Q \) is \( G_{p} \)-simple, we say that 
the associated hyperbolic fixed point
\( \alpha_{Q} \) is a 
\( G_{p} \)-\emph{simple}
(or \emph{simple}) number.
A hyperbolic number \( \alpha \)
is simple
if and only if
\( \alpha > 0 > \alpha^{\prime} \).
If 
\( \mathcal{A} \) is a hyperbolic equivalence class of 
\( \lambda \)-BQFs we write
\( \mathcal{Z_{A}} 
= \left\{\alpha_{Q} \mid Q \in \mathcal{A}, G_{p}\text{-simple}\right\} \),
the set of associated simple numbers.

Suppose that 
\( \alpha \) 
is a hyperbolic number and that 
\( \mathcal{A} \)
is a hyperbolic equivalence class of 
\( \lambda \)-BQFs.
Then
\( -\mathcal{A} 
	= \left\{ \left[ -A,-B,-C \right] 
	\mid \left[ A,B,C \right] \in \mathcal{A} \right\} \)
is also an equivalence class of 
\( \lambda \)-BQFs,
not necessarily distinct from
\( \mathcal{A} \).
A calculation shows that
\( Q_{\alpha^{\prime}} = -Q_{\alpha} \),
so
\( Q_{\alpha} \in \mathcal{A} \)
if and only if
\( Q_{\alpha ^{\prime}} \in \mathcal{-A} \).

\section{Rational period functions}

Choie and Zagier \cite{CZ93}
and Parson \cite{Par93}
gave an explicit characterization of RPFs on 
the modular group that 
made possible the correspondence in \cite{CR96}.
The problem of finding a characterization of 
RPFs on all of 
the Hecke groups
remains open,
although 
several authors have
done work in this direction
\cite{CR01, MR84, Res08, Sch93, SS95}.
These results provide enough information about
the structure of RPFs on the Hecke groups
for us to prove our correspondence theorem.

In this section we summarize the results 
about RPFs on the Hecke groups.
We give properties of RPF poles, 
quote a characterization of RPFs for the case we are considering,
and describe modifications that emphasize the second relation
\eqref{rpfrelation2}.

Throughout this section we suppose that \(q\) is an RPF of weight 
\(2k \in 2\Z^{+} \) on \(G_{p}\)
with pole set \(\Pset = \Pset(q) \).
We will not assume that \( k \) is odd or that 
\( q \) has a Hecke-symmetric set of poles until
subsection \ref{subsec:RPFchar};
the results of subsection \ref{subsec:poles} hold for all
RPFs
on \( G_{p} \).

\subsection{Poles}
\label{subsec:poles}

Hawkins \cite{Haw} introduced the idea of an 
\emph{irreducible system of poles}
(or
\emph{irreducible pole set}),
the minimal set of RPF poles, which must occur together because of the 
relations \eqref{rpfrelation1} and \eqref{rpfrelation2}.
Meier and Rosenberger
\cite{MR84}
observed that all the poles of
\( q \) are real and
Schmidt
\cite{Sch93}
proved that all nonzero poles are hyperbolic fixed points of 
\( G_{p} \).

Let 
\( \Pset^{+} \) and \( \Pset^{-} \)
denote the positive and negative poles, respectively,
in \( \Pset \),
and put 
\( \Pset^{*} = \Pset^{+} \cup \Pset^{-} \).
Schmidt 
\cite{Sch93} showed that 
\( \Pset^{*} \)
is a disjoint union of ``cycle pairs.''
We show in 
\cite{CR01}
that each of Schmidt's cycle pairs can be written as
\( {\mathcal{Z_{A}}} \cup T{\mathcal{Z_{A}}} \),
where \( \mathcal{A} \) is a hyperbolic equivalence class of 
\( \lambda \)-BQFs and
\( T\mathcal{Z_{A}} 
= \left\{ T\alpha \mid \alpha \in \mathcal{Z_{\mathcal{A}}} \right\} \).
Each
\( {\mathcal{Z_{A}}} \cup T{\mathcal{Z_{A}}} \)
is an irreducible system of poles;
the poles in 
\( \mathcal{Z_{A}} \) are positive and the poles in 
\( T{\mathcal{Z_{A}}} \) are negative.
In \cite{CR01} we also rewrite the negative poles in an irreducible system of poles as
\( T\mathcal{Z_{A}} = \mathcal{Z^{\prime}_{-A}} \),
where 
\( \mathcal{Z^{\prime}_{A}} 
= \left\{ \alpha^{\prime} \mid \alpha \in \mathcal{Z_{\mathcal{A}}} \right\} \).
As a result we may write any irreducible system of poles as
\( \mathcal{Z_{A}} \cup \mathcal{Z^{\prime}_{-A}} \)
for some hyperbolic \( \lambda \)-BQF equivalence class 
\( \mathcal{A} \).
A set
\( \mathcal{Z_{A}} \cup \mathcal{Z^{\prime}_{-A}} \)
has Hecke symmetry if and only if 
\( \mathcal{A} = -\mathcal{A} \), so
a Hecke symmetric irreducible systems of poles has the form
\( \mathcal{Z_{A}} \cup \mathcal{Z_{A}^{\prime}} \).
If all of the irreducible systems of poles of a RPF  are Hecke-symmetric we say that the RPF
itself is Hecke-symmetric.

\subsection{RPF characterization}
\label{subsec:RPFchar}

Meier and Rosenberger \cite{MR84} showed that
if an RPF of weight \( 2k \in 2\Z^{+} \) on \( G_{p} \) 
has a pole only at zero,
it must be of the form
\begin{equation}
  q_{0}(z) = 
    \label{RPFpoleatzero}
    \begin{cases}
      \nu(1-z^{-2k}), & \text{if } 2k \neq 2, \\
      \nu(1-z^{-2}) + \eta z^{-1}, & \text{if } 2k = 2,
    \end{cases}
\end{equation}
for constants \( \nu \) and \( \eta \).

Let \( k \) be an odd positive integer
and write 
\( Q_{\alpha}(z) = Q_{\alpha}(z,1) \)
for each \( \lambda \)-BQF \( Q \).
Suppose that \( q \) is a Hecke-symmetric RPF of weight 
\( 2k \) on \( G_{p} \).
In \cite{CR01} we show that
every such RPF
is of the form
\begin{equation}
    \label{}
    q(z) = \sum_{\ell=1}^{L} d_{\ell}
           \sum_{\alpha \in \mathcal{Z}_{\mathcal{A}_{\ell}}} 
       Q_{\alpha}(z)^{-k}
     + c_{0}q_{0}(z),
     \nonumber
\end{equation}
where 
each
\( \mathcal{A}_{\ell} \) 
is a hyperbolic
\( G_{p} \)-equivalence class of 
\( \lambda \)-BQFs,
the \( d_{\ell} \) \( (1 \leq \ell \leq L) \) are constants,
and \( q_{0}(z) \) is given by \eqref{RPFpoleatzero}.
We will see that the part of \( q \) with nonzero poles determines the
remainder term, 
so we write 
\begin{equation}
    \label{RPFsymmetricparts}
    q(z) = q^{*}(z) + c_{0}q_{0}(z),
\end{equation}
where
\begin{equation}
    \label{RPFsymmetricodd1*}
    q^{*}(z) = \sum_{\ell=1}^{L} d_{\ell}
           \sum_{\alpha \in \mathcal{Z}_{\mathcal{A}_{\ell}}} 
       Q_{\alpha}(z)^{-k}
\end{equation}
has poles
\begin{equation*}
    P^{*} = \bigcup_{\ell=1}^{L}
            \left( \mathcal{Z}_{\mathcal{A}_{\ell}} \cup 
                   \mathcal{Z}_{\mathcal{A}_{\ell}}^{\prime} \right).
\end{equation*}
Since both 
\( q \) and \( c_{0}q_{0} \) are RPFs, 
\( q^{*} \) is an RPF as well
and satisfies the relations
\eqref{rpfrelation1} and \eqref{rpfrelation2}.
We show in \cite[Lemma 9]{CR01} that 
\begin{equation}
    \label{RPFterm}
    \frac{D^{k/2}}{Q_{\alpha}(z)^{k}} 
    = \frac{(\alpha - \alpha^{\prime})^{k}}
      {(z-\alpha)^{k}(z-\alpha^{\prime})^{k}},
\end{equation}
where
\( D \) is the discriminant of the \( \lambda \)-BQF
\( Q_{\alpha} \).
Thus we may write 
\( q^{*} \) and \( q \) more explicitly as
\begin{equation}
    \label{RPFsymmetricoddexplicit*}
    q^{*}(z) = \sum_{\ell=1}^{L} d_{\ell} D_{\ell}^{-k/2}
           \sum_{\alpha \in \mathcal{Z}_{\mathcal{A}_{\ell}}} 
	       \frac{(\alpha - \alpha^{\prime})^{k}}
	      {(z-\alpha)^{k}(z-\alpha^{\prime})^{k}},
\end{equation}
and
\begin{equation}
    \label{RPFsymmetricoddexplicit}
    q(z) = \sum_{\ell=1}^{L} d_{\ell} D_{\ell}^{-k/2}
           \sum_{\alpha \in \mathcal{Z}_{\mathcal{A}_{\ell}}} 
	       \frac{(\alpha - \alpha^{\prime})^{k}}
	      {(z-\alpha)^{k}(z-\alpha^{\prime})^{k}}
	      + c_{0}q_{0}(z).
\end{equation}

\section{The direct theorem}

In this section we show that an entire automorphic integral of positive even weight on one of the Hecke groups may be associated with a Dirichlet series satisfying a functional equation.
We restrict our attention to automorphic integrals of weight 
\( 2k \)
(\( k \) odd)
with Hecke-symmetric RPFs.
We show that the Dirichlet series
functional equation involves a remainder term, which comes from the RPF for the automorphic integral.

Let \( k \) be an odd positive integer.
Suppose that
$F$ is an 
entire automorphic integral of weight 
\( 2k \in 2\Z^{+} \) 
on \( G_{p} \)
with Hecke-symmetric RPF \( q \).
We may assume without loss of generality that 
\( F \) is a 
\emph{cusp} automorphic integral, that is,
\( a_{0} = 0 \) in the Fourier expansion
\eqref{Fourierexp}.

Write $z=x+iy$ with $x,y \in \R$.
It can be shown
\cite[622-623]{Kno74}
that $F$ satisfies
\begin{equation}
	\label{MIgrowth1}
	\left| F(z) \right| \leq K \left(|z| ^{\alpha}
			+ y^{- \beta} \right),
			\; \;   z \in \mathcal{H}
\end{equation}
for some positive real numbers $K$, $\alpha$ and $\beta$.
It follows that the coefficients $a_{n}$ in the Fourier expansion (\ref{Fourierexp}) for $F$ satisfy
\begin{equation}
	a_{n} = {\cal O} (n^{\beta}), \; \; n \rightarrow +\infty.
	\label{coeffgrowth}
\end{equation}
This, with $a_{0} = 0$ in (\ref{Fourierexp}), implies that
\begin{equation}
	\label{MIgrowth2}
	F(iy) = {\cal O} (e^{-2 \pi y/\lambda}), \; \; y \rightarrow +\infty.
\end{equation}
The growth estimates
(\ref{MIgrowth1}) and (\ref{MIgrowth2}) allow us to
define the
\emph{Mellin transform of $F$},
\begin{equation}
	\label{PhiIntegral}
	\Phi (s) = \int_{0}^{\infty} F(iy)y^{s}\frac{dy}{y},
\end{equation}
a function of the complex variable $s = \sigma + it$.
This integral converges for 
$\sigma > \beta$.
For $\sigma > \beta +1$, we can integrate term by term to get
\begin{equation}
	\label{PhiDS}
	\Phi (s) = \left( \frac{2\pi}{\lambda} \right)^{-s} \Gamma (s) \phi (s),
\end{equation}
where
\begin{equation}
	\label{phi}
	\phi (s) = \sum_{n=1}^{\infty} \frac{a_{n}}{n^{s}}
\end{equation}
is the {\em Dirichlet series associated with $F$\/}.
The bound on the growth of the Fourier coefficients 
$a_{n}$ (\ref{coeffgrowth}) 
implies that the sum in (\ref{phi}) converges absolutely and uniformly on compact subsets of the right half plane $\sigma > \beta + 1$,
so that $\phi(s)$ is analytic there.

We invert part of the Mellin transform of \( F \) and use 
the automorphic relation \eqref{mirelation}
and the RPF decomposition \eqref{RPFsymmetricparts}
to get
\begin{equation*}
\begin{split}
	\int_{0}^{1} F(iy)y^{s}\frac{dy}{y} 
	& = 
	\int_{1}^{\infty} F\left( \frac{-1}{iy} \right)y^{-s}\frac{dy}{y} \\
	& = 
	-\int_{1}^{\infty} F(iy)y^{2k-s}\frac{dy}{y}
	-\int_{1}^{\infty} c_{0}q_{0}(iy)y^{2k-s}\frac{dy}{y} \\
	& \phantom{xxxxxxxxxxxxxxxxxxxxxxxxxxxxx} 
	-\int_{1}^{\infty} q^{*}(iy)y^{2k-s}\frac{dy}{y}.
\end{split}
\end{equation*}
Thus
\begin{equation}
	\label{}
	\Phi(s) = D(s) + E^{0}(s) + E^{*}(s),
	\nonumber
\end{equation}
where
\begin{equation}
	\label{Phientire}
	D(s) = \int_{1}^{\infty} F(iy)\left[ y^{s} -y^{2k-s} \right]\frac{dy}{y},
\end{equation}
\begin{equation}
    \label{PhiRPF0}
    E^{0}(s) = -\int_{1}^{\infty} c_{0}q_{0}(iy)y^{2k-s}\frac{dy}{y},
\end{equation}
and
\begin{equation}
    \label{PhiRPF*}
    E^{*}(s) = -\int_{1}^{\infty} q^{*}(iy)y^{2k-s}\frac{dy}{y}.
\end{equation}
Now 
\( D(s) \) is entire
and satisfies the functional equation
\begin{equation}
    \label{PhientireFE}
    D(2k-s)+D(s)=0.
\end{equation}
The expression
\eqref{RPFpoleatzero} 
for \( q_{0} \)
implies that 
\( q_{0}(z) = \mathcal{O}(1) \) as
\( |z| \rightarrow \infty \),
so the integral defining 
\( E^{0}(s) \) in 
\eqref{PhiRPF0} converges in the right half-plane 
\( \sigma > 2k \).
The expression
\eqref{RPFsymmetricodd1*} 
for \( q^{*} \)
implies that 
\( q^{*}(z) = \mathcal{O}\left( |z|^{-2k} \right) \) as
\( |z| \rightarrow \infty \),
so the integral defining 
\( E^{*}(s) \) in 
\eqref{PhiRPF*} converges in the right half-plane 
\( \sigma > 0 \).

In order to write the functional equation for
\( \Phi(s) \) 
that is suggested by
\eqref{PhientireFE} we need meromorphic continuations of 
\( E^{0}(s) \) and \( E^{*}(s) \) to the \( s \)-plane.
We use \eqref{RPFpoleatzero} to calculate that
\begin{equation}
  E^{0}(s) = 
    \label{PhiRPF0explicit}
    \begin{cases}
      -\tilde{a}_{0}\left( \frac{1}{s-2k} + \frac{1}{s} \right),
	       & \text{if } 2k \neq 2, \\
      -\tilde{a}_{0}\left( \frac{1}{s-2} + \frac{1}{s} \right)
		+ \frac{\tilde{b}_{1}i}{s-1}, 
			& \text{if } 2k = 2,
    \end{cases}
\end{equation}
where
\( \tilde{a}_{0} = a_{0}c_{0} \)
and
\( \tilde{b}_{0} = b_{0}c_{0} \).
In every case \( E^{0}(s) \)
has a meromorphic continuation to all of 
the complex \( s \)-plane,
with simple poles at
\( s = 0, 2k \) 
(and at \( s=1 \) if \( 2k=2 \)).
Furthermore, 
\( E^{0}(s) \) satisfies the same function equation
as \( D(s) \),
\begin{equation}
    \label{}
    E^{0}(2k-s) + E^{0}(s) = 0.
    \nonumber
\end{equation}

For the meromorphic continuation of
\( E^{*}(s) \)
we first use a partial fraction decomposition 
of the right side of \eqref{RPFterm}
\begin{equation*}
    \frac{(\alpha - \alpha^{\prime})^{k}}{(z-\alpha)^{k}(z-\alpha^{\prime})^{k}}
    = \sum_{m=1}^{k} \frac{a_{m,\alpha}}{(z-\alpha)^{m}}
   + \sum_{n=1}^{k} \frac{b_{n,\alpha^{\prime}}}{(z-\alpha^{\prime})^{n}},
\end{equation*}
where 
\( a_{m,\alpha} = (-1)^{m-k} \binom{2k-m-1}{k-1} 
	(\alpha - \alpha^{\prime})^{m-k} \)
and 
\( b_{n,\alpha^{\prime}} = (-1)^{k} \binom{2k-n-1}{k-1} 
	\left( \alpha - \alpha^{\prime} \right)^{n-k} \).
This along with
\eqref{RPFsymmetricodd1*} gives us
\begin{equation}
    \label{RPFsymmetricodd2*}
    q^{*}(z) = \sum_{\ell=1}^{L} c_{\ell}
           \sum_{\alpha \in \mathcal{Z}_{\mathcal{A}_{\ell}}} 
	       \left( \sum_{m=1}^{k} \frac{a_{m,\alpha}}{(z-\alpha)^{m}}
	   + \sum_{n=1}^{k} 
	   \frac{b_{n,\alpha^{\prime}}}{(z-\alpha^{\prime})^{n}} \right),
\end{equation}
where
\( c_{\ell} = d_{\ell}D_{\ell}^{-k/2} \)
with 
\( D_{\ell} \)
the discriminant of the BQFs in the equivalence class
\( \mathcal{A}_{\ell} \).

We use the integral representation for the hypergeometric function
\cite[equation (9.1.6)]{Leb72}
\begin{equation}
    \label{hypergeometricintegral}
    _{2}F_{1}[a,b;c;z] 
	    = \frac{\Gamma(c)}{\Gamma(b)\Gamma(c-b)}
	    \int_{0}^{1} y^{b-1}(1-y)^{c-b-1}(1-yz)^{-a} dy,
\end{equation}
for 
\( \re(c) > \re(b) >0 \) 
and
\( \left| \arg (1-z) \right| < \pi \).
We let
\( a=m \), 
\( b=m-s \), and
\( c=1+m-s \)
for \( s \) a complex variable and
\( m \) a positive integer,
and we use a
change of variables to invert, so
\begin{equation*}
    \int_{1}^{\infty} \frac{y^{s}}{(y-z)^{m}} \frac{dy}{y}
    = \frac{1}{r-s} {} _{2}F_{1}[m,m-s;1+m-s;z],
\end{equation*}
for 
\( \sigma < m \) 
and
\( \left| \arg (1-z) \right| < \pi \).
We let 
\( z = i\alpha \) for \( \alpha \in \R \)
and multiply by \( i^{-m} \), and get
\begin{equation}
    \label{integralformula1}
    \int_{1}^{\infty} \frac{y^{s}}{(iy-\alpha)^{m}} \frac{dy}{y}
    = \frac{i^{-m}}{m-s} {} _{2}F_{1}[m,m-s;1+m-s;-i\alpha],
\end{equation}
for 
\( \sigma < m \).

The expression \eqref{RPFsymmetricodd2*}
and formula \eqref{integralformula1} together imply that
\( E^{*}(s) \) is a finite linear combination of terms of the form
\begin{equation}
    \label{integralformula2}
    \int_{1}^{\infty} \frac{y^{2k-s}}{(iy-\alpha)^{m}}
    \frac{dy}{y}
    = \frac{i^{-m}}{s-2k+m} {} 
    _{2}F_{1}[m,s-2k+m;s-2k+m+1;-i\alpha],
\end{equation}
for \( \alpha \in \Pset^{*} \),
\( 1 \leq m \leq k \),
and 
\( \sigma > 2k-m \).
Now the hypergeometric function
\( _{2}F_{1}[a,b;c;z] \) is an entire function of
\( a \) and \( b \), and
a meromorphic function of \( c \)
with simple poles at 
\( c = 0, -1, -2, \dots \).
Thus the function in 
\eqref{integralformula2}
is meromorphic in the \( s \)-plane
with simple poles at
\( s = 2k-m, 2k-m-1, 2k-m-2, \dots \).
Since  
\( E^{*}(s) \) is a finite linear combination of terms of the form
\eqref{integralformula2} with
\( 1 \leq m \leq k \),
we have that 
\( E^{*}(s) \) is a meromorphic function
with simple poles at most at
\( s = 2k-1, 2k-2, \dots \).
Thus \( \Phi(s) \) is meromorphic in the
\( s \)-plane.

We now show that 
\( \Phi(s) \)
is bounded in lacunary vertical strips
of the form
\begin{equation}
    \label{lacunaryVS}
    S(\sigma_{1},\sigma_{2};t_{0}) 
    = \left\{ s = \sigma + it \mid 
    \sigma_{1} \leq \sigma \leq \sigma_{2}, |t| \geq t_{0} >0 \right\}.
\end{equation}
Now 
\( D(s) \) is bounded in every vertical strip
since the integral in 
\eqref{Phientire} converges for every \( s \),
\( |y^{s}| = y^{\sigma} \geq y^{\sigma_{2}} \)
and
\( |y^{2k-s}| = y^{2k-\sigma} \leq y^{2k-\sigma_{1}} \),
and the integrand is unchanged as 
\( t \rightarrow \infty \).
Also,
the expression \eqref{PhiRPF0explicit} shows that
\( E^{0}(s) \) is bounded in every
lacunary vertical strip,
since its poles 
are excluded from every 
\( S(\sigma_{1},\sigma_{2};t_{0}) \) and
each term 
approaches zero
as 
\( t \rightarrow \infty \).

It remains to show that 
\( E^{*}(s) \) is bounded in every
lacunary vertical strip.
Since
\( E^{*}(s) \) is a finite sum of terms of the form
\eqref{integralformula2},
it is sufficient to 
prove that for every purely imaginary
\( \beta \),
and for \( 1 \leq m \leq k \)
the function
\begin{equation*}
    \frac{1}{s-2k+m} {} 
    _{2}F_{1}[m,s-2k+m;s-2k+m+1;\beta]
\end{equation*}
is bounded in any \( S(\sigma_{1},\sigma_{2};t_{0}) \).
Given any lacunary vertical strip of the form 
\eqref{lacunaryVS}
we use 
\eqref{hypergeometricintegral}
to write
\begin{multline*}
    \frac{1}{s-2k+m} {} _{2}F_{1}[m,s-2k+m;s-2k+m+1;\beta] \\
	    = \int_{0}^{1} y^{s-2k+m-1}(1-\beta y)^{-m} dy,
\end{multline*}
for 
\( \sigma > 2k-m \) and
\( \alpha \in \R \).
This integral is bounded in every
\( S(\sigma_{1},\sigma_{2};t_{0}) \)
for which \( \sigma_{1} > 2k-m \), since
\( \left| y^{s-2k+m-1}(1-\beta y)^{-m} \right| \leq y^{\sigma_{1}-2k+m-1} \)
for \( s \geq \sigma_{1} \)
and \( 0 \leq y \leq 1 \).
If \( \sigma_{1} \leq 2k-m \) we let 
\( n \in \Z^{+} \) such that 
\( \sigma_{1} > 2k-m-n \)
and we integrate by parts \( n \) times.
The result is
\begin{equation}
\begin{split}
    \label{PhiRPF*term3}
	\frac{1}{s-2k+m}& {} _{2}F_{1}[m,s-2k+m;s-2k+m+1;\beta] \\
	& = \sum_{j=1}^{n} \frac{(-1)^{j+1} \Gamma(s-2k+m) m! \beta^{j-1}}
		    {\Gamma(s-2k+m) (1-\beta)^{m+j-1} (m-j+1)!}  \\
	& \phantom{xx} 
		+ \frac{(-1)^{n+1} \Gamma(s-2k+m) m! \beta^{n}}
		    {\Gamma(s-2k+m+n) (m-n)!}
	    \int_{0}^{1} y^{s-2k+m+n-1}(1-\beta y)^{-m-n} dy,
\end{split}
\end{equation}
for \( \sigma > 2k-m-n \).
This integral is bounded in
\( S(\sigma_{1},\sigma_{2};t_{0}) \) since
\( \left| y^{s-2k+m+n-1}(1-\beta y)^{-m-n} \right| \leq y^{\sigma_{1}-2k+m+n-1} \)
for \( s \geq \sigma_{1} \)
and \( 0 \leq y \leq 1 \).
The other expressions 
on the right-hand side of 
\eqref{PhiRPF*term3} are rational in 
\( s \) with simple poles at integer values, 
and they each approach zero as 
\( t \rightarrow \infty \).
Thus the function in 
\eqref{PhiRPF*term3}
is bounded in every
\( S(\sigma_{1},\sigma_{2};t_{0}) \),
which implies that 
\( E^{*}(s) \) and \( \Phi(s) \) are bounded there as well.

Since
\( \Phi(s) \)
is meromorphic
we may write the functional equation suggested by 
\eqref{PhientireFE},
\begin{equation}
    \label{PhiFE}
    \Phi(2k-s) + \Phi(s) = R(s),
\end{equation}
where
\( R(s) \)
is a meromorphic function we call the
\emph{remainder term}.
By
\eqref{PhientireFE} and \eqref{EzeroFE}
we have
\begin{equation}
    \label{remainderterm}
    R(s) = E^{*}(2k-s) + E^{*}(s),
\end{equation}
so the remainder term depends only on 
\( q^{*}(z) \),
the part of the RPF with nonzero poles.
The expression 
\eqref{PhiFE}
(or \eqref{remainderterm})
implies that \( R(s) \) satisfies the
(first) relation 
\begin{equation}
    \label{Rrelation1}
    R(2k-s) - R(s) = 0.
\end{equation}

We must calculate an explicit expression for
\( R(s) \), 
in order to give meaning to
\eqref{PhiFE}
and to prove the converse theorem.
If we use the fact that \( q^{*} \) satisfies the first relation
\eqref{rpfrelation1} to replace 
\( q^{*}(iy) \) in \eqref{PhiRPF*} and invert, 
we have
\begin{equation}
    \label{PhiRPF*alt1}
    E^{*}(s) = - \int_{0}^{1} q^{*}(iy)y^{s}\frac{dy}{y}.
\end{equation}
On the other hand, if we substitute directly into
\eqref{PhiRPF*} we have
\begin{equation}
    \label{PhiRPF*alt2}
    E^{*}(2k-s) = -\int_{1}^{\infty} q^{*}(iy)y^{s}\frac{dy}{y}.
\end{equation}
The integral in \eqref{PhiRPF*alt1}
converges for
\( \sigma > 0 \) since 
\( q^{*}(iy) \) is bounded as
\( y \rightarrow 0 \), and the
integral in \eqref{PhiRPF*alt2}
converges for
\( \sigma < 2k \) since 
\( q^{*}(iy) = \mathcal{O}(y^{-2k}) \) as
\( y \rightarrow \infty \).
Thus for \( 0 < \sigma < 2k \) we have
\begin{equation}
    \label{remainderterm2}
    R(s) = -\int_{0}^{\infty} q^{*}(iy)y^{s}\frac{dy}{y},
\end{equation}
that is,
\( R(s) \) is the negative of the Mellin transform
of \( q^{*}(z) \).
This expression makes it clear that the first relation
for \( q^{*}(z) \) 
leads directly to the first relation
for \( R(s) \).
If we use \eqref{rpfrelation1} to replace
\( q^{*}(iy) \) in \eqref{remainderterm2} and invert the variable of integration
we get \eqref{Rrelation1}.

We will substitute the expression for 
\( q^{*} \) given by \eqref{RPFsymmetricoddexplicit*} into 
\eqref{remainderterm2} and write
\( R(s) \) as a linear combination of integrals of the form
\begin{equation}
    \label{}
    \int_{0}^{\infty} \frac{y^{s}}
      {\left( iy - \alpha \right)^{k} \left( iy - \alpha^{\prime} \right)^{k}}
      \frac{dy}{y},
      \nonumber
\end{equation}
which converge for \( 0 < \sigma < 2k \).
The evaluation of these integrals involves exponential functions of the form
\( z^{a} = e^{a \log z} \),
where
$\log z = \log |z| + i\arg z$
for
$z \in {\bf C}$.
We will take the principal branch for each logarithm, using the convention that
$-\pi \leq \arg z < \pi$.
We need the integral formula in the following lemma,
which uses the beta function 
\( B(a,b) = \frac{\Gamma(a)\Gamma(b)}{\Gamma(a+b)} \) and 
the hypergeometric function.

\begin{lemma}
\label{lemma:MTformulasym}
	Let \( k \in \Z^{+} \), and let  
	$\delta$ and $\epsilon$ be nonzero real numbers that satisfy one of:
	$\delta<0<\epsilon$, $\epsilon<0<\delta$, $0<\epsilon<\delta$, or
	$\delta<\epsilon<0$.
	Then
	\begin{eqnarray}
		\lefteqn{
		\int_{0}^{\infty}
			\frac{y^{s}}{(iy-\delta)^{k}(iy-\epsilon)^{k}} 
		\frac{dy}{y}
		}	\nonumber		\\
	& = &
		i^{s} \delta^{s-k} B(2k-s,s-k) (\delta - \epsilon)^{-k}
		{} _{2}F_{1}\left[ k, 1-k; k-s+1;\frac{\epsilon}{\epsilon-\delta} \right]
		\nonumber		\\
	&  & 
	+	i^{s} \epsilon^{s-k} B(s,k-s) (\delta - \epsilon)^{-k}
		{} _{2}F_{1}\left[ k, 1-k; s-k+1;\frac{\epsilon}{\epsilon-\delta} \right],
	\label{MTformulasym}
	\end{eqnarray}
	for $0<\mbox{Re} \, s<2k$.
\end{lemma}

\begin{proof}
At several places the proof below involves branching questions
that reduce to calculations of the form
\( \left( z_{1}z_{2} \right)^{s} = z_{1}^{s}z_{2}^{s} \),
which is valid if
\( \arg(z_{1}) + \arg(z_{2}) = \arg\left( z_{1}z_{2} \right) \)
using our argument convention.

We start with \eqref{hypergeometricintegral} and
change variables by letting 
\( y = \frac{u}{u+1} \).
We also let 
\( z = 1-v/w \), where 
$v$ is a positive real number and 
$w \in {\bf C} \setminus {\bf R}$,
which implies that
\( |\arg(1-z)| = \left| \arg \left( v/w \right) \right| < \pi \).
The result is
\begin{displaymath}
	\int_{0}^{\infty} u^{b-1} (u+1)^{a-c} (w+vu)^{-a} du
=
	w^{-a} B(b,c-b)
	\hspace{.5em} F_{\hspace{-1.03em} 2 \hspace{.65em} 1}
	[a,b;c;1-v/w],
\end{displaymath}
for
$\mbox{Re} \, c > \mbox{Re} \, b > 0$, 
$v > 0$ and 
$w \in {\bf C} \setminus {\bf R}$.
We change variables again by letting 
$u=y/v$, 
and put $a=k$, $b=s$, and $c=2k$
($k \in {\bf Z}^{+}$ and $s \in {\bf C}$)
so that 
\begin{equation*}
	\int_{0}^{\infty} 
	\frac{y^{s}}{(y+v)^{k}(y+w)^{k}}
	\frac{dy}{y}
=
	v^{s-k} w^{-k} B(s,2k-s)
	\hspace{.5em} F_{\hspace{-1.03em} 2 \hspace{.65em} 1}
	[k,s;2k;1-v/w]
\end{equation*}
for $0< \sigma = \mbox{Re} \, (s) < 2k$,
$v > 0$ and 
$w \in {\bf C} \setminus {\bf R}$.
We will apply the two identities for hypergeometric functions 
\cite[equations (9.5.3) and (9.5.9)]{Leb72},
\begin{equation}
	\hspace{.5em} F_{\hspace{-1.03em} 2 \hspace{.65em} 1}
	[a,b;c;z]
=
	(1-z)^{c-a-b}
	\hspace{.5em} F_{\hspace{-1.03em} 2 \hspace{.65em} 1}
	[c-a,c-b;c;z],
\label{2F1:3}
\end{equation}
for $|\arg(1-z)|<\pi$, and
\begin{eqnarray}
	\hspace{.5em} F_{\hspace{-1.03em} 2 \hspace{.65em} 1}
	[a,b;c;z]
& = &
	(-z)^{-a} \frac{\Gamma(c) \Gamma(b-a)}{\Gamma(c-a) \Gamma(b)}
	\hspace{.5em} F_{\hspace{-1.03em} 2 \hspace{.65em} 1}
	[a,1+a-c;1+a-b;1/z]
	\nonumber		\\
&  &	+
	(-z)^{-b} \frac{\Gamma(c) \Gamma(a-b)}{\Gamma(c-b) \Gamma(a)}
	\hspace{.5em} F_{\hspace{-1.03em} 2 \hspace{.65em} 1}
	[b,1+b-c;1+b-a;1/z],
	\nonumber		\\
&  &
\label{2F1:2}
\end{eqnarray}
for $|\arg(-z)|<\pi$ and $|\arg(1-z)|<\pi$.
We apply (\ref{2F1:2}) with 
$a=k$, $b=s$, $c=2k$ and 
$z = 1- v/w$,
and get
\begin{eqnarray}
\lefteqn{
	\int_{0}^{\infty} 
	\frac{y^{s}}{(y+v)^{k}(y+w)^{k}}
	\frac{dy}{y}
	}		\nonumber		\\
& = &
	v^{s-k} (v-w)^{-k} B(2k-s,s-k)
	\hspace{.5em} F_{\hspace{-1.03em} 2 \hspace{.65em} 1}
	\left[k,1-k;k-s+1;\frac{w}{w-v}\right]
	\nonumber		\\
&  & +
	v^{s-k} w^{s-k} (v-w)^{-s} B(s,k-s)
	\hspace{.5em} F_{\hspace{-1.03em} 2 \hspace{.65em} 1}
	\left[s,1+s-2k;s-k+1;\frac{w}{w-v}\right],
	\nonumber	
\end{eqnarray}
for $0< \sigma < 2k$,
$v > 0$ and 
$w \in {\bf C} \setminus {\bf R}$.
Next we apply (\ref{2F1:3}) to the hypergeometric function in the second term, with
$a=s$, $b=1+s-2k$, $c=s-k+1$ and
$z = \frac{w}{w-v}$.
After simplifying we have
\begin{eqnarray}
\lefteqn{
	\int_{0}^{\infty} 
	\frac{y^{s}}{(y+v)^{k}(y+w)^{k}}
	\frac{dy}{y}
	}		\nonumber		\\
& = &
	v^{s-k} (v-w)^{-k} B(2k-s,s-k)
	\hspace{.5em} F_{\hspace{-1.03em} 2 \hspace{.65em} 1}
	\left[k,1-k;k-s+1;\frac{w}{w-v}\right]
	\nonumber		\\
&  &
	+ w^{s-k} (v-w)^{-k} B(s,k-s)
	\hspace{.5em} F_{\hspace{-1.03em} 2 \hspace{.65em} 1}
	\left[k,1-k;s-k+1;\frac{w}{w-v}\right],
\label{2F1:4}
\end{eqnarray}
for $0<\sigma<2k$, 
$v>0$ and 
$w \in {\bf C} \setminus {\bf R}$.
Next we restrict $w$ to 
$\arg w = \pm \pi/2$
and put
$v = \delta i$
and 
$w = \epsilon i$
in (\ref{2F1:4}).
A simplification gives us
(\ref{MTformulasym}) for 
$\arg \delta = -\pi/2$
and
$\epsilon \in {\bf R}$, $\epsilon \neq 0$.

To complete the proof we must consider the possible values of $\delta$ and $\epsilon$.
We consider $\delta$ to be a complex variable and do an analytic continuation of (\ref{MTformulasym}) in $\delta$ to the region 
$-\pi \leq \arg \delta < \pi/2$ subject to the restriction that 
$\left| \arg \left( 1 - \left( \frac{\epsilon}
	{\epsilon-\delta} \right) \right) \right|
=	\left| \arg \left( \frac{\delta}
	{\delta-\epsilon} \right) \right| < \pi$.
If we let $\delta \in {\bf R}$, this restriction is
$\frac{\delta}{\delta-\epsilon} > 0$,
which implies the restrictions given in the statement of the Lemma.
\end{proof}

We calculate 
\( R(s) \) explicitly by
substituting \eqref{RPFsymmetricoddexplicit*} into 
\eqref{remainderterm2} to get
\begin{equation*}
    R(s) = -\sum_{\ell=1}^{L} d_{\ell} D_{\ell}^{-k/2}
           \sum_{\alpha \in \mathcal{Z}_{\mathcal{A}_{\ell}}} 
               \left( \alpha - \alpha^{\prime} \right)^{k}
               \int_{0}^{\infty}
               \frac{y^{s}}{(iy-\alpha)^{k}(iy-\alpha^{\prime})^{k}}
               \frac{dy}{y},
\end{equation*}
for \( 0 < \sigma < 2k \)
with 
\( D_{\ell} \) the discriminant of the 
\( \lambda \)-BQFs in 
\( \mathcal{A}_{\ell} \).
We may apply Lemma \ref{lemma:MTformulasym}
to the integrals in this expression,
since 
\( \alpha > 0 > \alpha^{\prime} \).
This gives us
\begin{multline}
    \label{remainderterm4}
    R(s) = - \sum_{\ell=1}^{L} d_{\ell} D_{\ell}^{-k/2} i^{s}  \\
           \times \sum_{\alpha \in Z_{\mathcal{A}_{\ell}}} 
           \lbrace \alpha^{s-k}B(2k-s,s-k)
           _{2}F_{1}\left[ k,1-k;k-s+1;
           \frac{\alpha^{\prime}}{\alpha^{\prime}-\alpha} \right]  \\
           + \left( \alpha^{\prime} \right)^{s-k}B(s,k-s)
           _{2}F_{1}\left[ k,1-k;s-k+1;
           \frac{\alpha^{\prime}}{\alpha^{\prime}-\alpha} \right] \rbrace.
\end{multline}

We may use this expression
to verify directly that 
the remainder term
satisfies the first relation
\eqref{Rrelation1}.
We use \eqref{remainderterm4} to write
\( R(2k-s) \), then
use the fact that
\( \alpha \in  \mathcal{Z}_{\mathcal{A}}\)
implies \( -\frac{1}{\alpha} = \gamma^{\prime} \)
and 
\( -\frac{1}{\alpha^{\prime}} = \gamma \in  \mathcal{Z}_{\mathcal{A}} \),
as well as a calculation that 
\( \frac{\alpha^{\prime}}{\alpha^{\prime} - \alpha} 
=  \frac{\gamma^{\prime}}{\gamma^{\prime} - \gamma} \)
to obtain \( R(s) \).

We have proved the following theorem.
\begin{theorem}
	\label{directtheorem}
	Fix \( p \geq 3 \),
	let \( \lambda = \lambda_{p} \), and
	let \( k \) be an odd positive integer.
	Suppose that
	\( F(z) \) is an 
	entire automorphic integral of weight 
	\( 2k \) 
	on \( G_{p} \)
	with Hecke-symmetric rational period function \( q(z) \)
	given by 
	\eqref{RPFsymmetricodd1*}.
	Suppose that \( F \) has the Fourier expansion
	\eqref{Fourierexp} with zero constant term,
	and that \( \Phi(s) \) is given by
	\eqref{PhiIntegral} for 
	\( \re(s) > \beta \), 
	for some positive \( \beta  \).
	
	Then for
	\( \re(s) > \beta +1 \)
	\( \Phi(s) \) is also given by 
	\eqref{PhiDS} and \eqref{phi}, and
	\newline
	(a)  
	\( \Phi(s) \) has a meromorphic continuation to the 
	\( s \)-plane with, at worst, simple poles at integer points
	\( m \leq 2k \);
	for \( \re(s) > \beta \) we have
	\begin{equation*}
	    \Phi(s) = D(s) + E^{0}(s) + E^{*}(s);
	\end{equation*}
	\( D(s) \) is given by 
	\eqref{Phientire} and is entire, while
	\( E^{0}(s) \) and \( E^{*}(s) \) 
	are given by 
	\eqref{PhiRPF0} and
	\eqref{PhiRPF*} (respectively)
	and have meromorphic continuations to the 
	\( s \)-plane;
	furthermore,
	\newline
	(b)  
	\( \Phi(s) \) is bounded in every lacunary vertical strip of the form
	\eqref{lacunaryVS}, and
	\newline
	(c)  
	\( \Phi(s) \) satisfies the functional equation
	\eqref{PhiFE} with 
	\( R(s) \) given by 
	\eqref{remainderterm4}.
\end{theorem}

\section{The second relation}

We have observed that the remainder term 
\( R(s) \) satisfies one relation 
\eqref{Rrelation1}.
In this section we show that
\( R(s) \) must satisfy a 
second relation, which follows from the fact that the corresponding
RPF
\( q^{*}(z) \) satisfies the second relation
\eqref{rpfrelation2}.
We first modify
\eqref{RPFsymmetricodd1*} in order to exhibit the fact that 
\( q^{*} \) satifies
\eqref{rpfrelation2}.

We let 
\( U = ST 
	= \left( \begin{smallmatrix} \lambda&-1\\ 1&0 \end{smallmatrix} \right)
\)
and focus on the effect of \( U \) on the poles of an RPF.
Schmidt proves in \cite{Sch93} that if
\( \alpha \in \Pset^{*} \), then 
\( U^{m}\alpha \in \Pset^{*} \) for exactly one
\( m \), \( 1 \leq m \leq p-1 \).
Thus we may separate the poles in
\( \Pset^{*} \) into pairs that are images of each other under some power of
\( U \).
For each pair of poles we will determine the power 
\( m \) from the sizes of the poles.

A calculation shows that 
\[ 0 = U^{p}(0) < U^{p-1}(0) < \cdots < U^{2}(0) < U(0) 
= \infty, \]
so we may partition the extended real line
into a disjoint union of \( p \) 
half-open intervals,
\begin{equation}
	\label{linedecomp}
	[-\infty,\infty) = [-\infty,U^{p}(0) \, ) \cup 
	[U^{p}(0),U^{p-1}(0) \, ) 
	\cup \cdots \cup [U^{2}(0),U(0) \, ).
\end{equation}
\( U \) maps each interval to the previous interval,
the first to the last, and left endpoints to left endpoints.
We denote the \( j \)th interval 
of \eqref{linedecomp} by
\( I_{j} \), that is
\( I_{j} =  [U^{p-j+2}(0),U^{p-j+1}(0) \, )\)
for \( 1 \leq j \leq p \).
We note that all negative poles are in 
\( I_{1} \) and each positive pole is in
one of the \( I_{j} \),
\( 2 \leq j \leq p \).

\begin{lemma}
	\label{lemma:mappingpoles}
	Fix \( p \geq 3 \) and
	\( 2 \leq j \leq p \).
	For every 
	\( G_{p} \)-equivalence class ${\cal A}$ of 
	\( \lambda_{p} \)-BQFs we have
	\( \left\{ \beta^{\prime} \mid 
		\beta \in \mathcal{Z_{-A}} \cap I_{p-j+2} \right\}
					=
		\left\{ U^{j-1}\alpha \mid 
		\alpha \in \mathcal{Z_{A}} \cap I_{j} \right\}
					\),
	where \( I_{j} \) is the \( j \)th interval of \eqref{linedecomp}
	and 
	\( \beta^{\prime} \) is the Hecke conjugate of 
	\( \beta \).
\end{lemma}

\begin{proof}
	We show containment in both directions. 
	First suppose that 
	\( \gamma \in \left\{ \beta^{\prime} \mid 
		\beta \in \mathcal{Z_{-A}} \cap I_{p-j+2} \right\} \), so
	\( \gamma = \beta^{\prime} \) for some 
	\( \beta \) in both \( \mathcal{Z_{-A}} \) and \( I_{p-j+2} \).
	Let \( \alpha = U^{p-j+1}\gamma \), so
	\( \gamma = U^{j-1}\alpha \).
	We need to show that \( \alpha \)
	is in both \( \mathcal{Z_{A}} \) and
	\( I_{j} \).
	Now \( \beta \) is simple, so
	\( \gamma = \beta^{\prime} \in I_{1} \), which implies that
	\( \alpha = U^{-(j-1)}\gamma \in I_{j} \).
	We also have
	\( \gamma^{\prime} = \beta \in I_{p-j+2} \), so
	\( \alpha^{\prime} = U^{p-j+1}\gamma^{\prime} \in I_{1} \).
	Thus \( \alpha^{\prime} < 0 < \alpha \),
	so \( \alpha \) is simple.
	Now \( \beta \in \mathcal{Z_{-A}} \) implies that
	\( Q_{\beta} \in \mathcal{-A} \), so
	\( Q_{\gamma} = Q_{\beta^{\prime}} \in \mathcal{A} \).
	But since \( \alpha = U^{p-j+1}\gamma \) this means that
	\( Q_{\alpha} \in \mathcal{A} \).
	Thus
	\( \alpha \in \mathcal{Z_{A}} \)
	and we have shown that
	\( \left\{ \beta^{\prime} \mid 
		\beta \in \mathcal{Z_{-A}} \cap I_{p-j+2} \right\}
					\subseteq
		\left\{ U^{j-1}\alpha \mid 
		\alpha \in \mathcal{Z_{A}} \cap I_{j} \right\} \).
	
	The demonstration of containment in the other direction is a similar calculation.
\end{proof}

Lemma \ref{lemma:mappingpoles} means that 
every negative pole of an RPF is connected by a power of 
\( U \) to a positive pole.
Specifically,
if
\( \beta^{\prime} \in \Pset^{-} \),
then 
\( \beta^{\prime} = U^{j-1}\alpha \) for some
\( \alpha \in \Pset^{+} \).
The power \( j-1 \) is determined by the size of 
\( \beta \), the Hecke conjugate of 
\( \beta^{\prime} \).
Although 
Lemma \ref{lemma:mappingpoles} holds for all equivalence classes of 
\( \lambda_{p} \)-BQFs, 
we will apply it to equivalence classes that have Hecke-symmetry.
In this setting equivalence classes satisfy
\( -\mathcal{A} = \mathcal{A} \), 
so we have
	\( \left\{ \beta^{\prime} \mid 
		\beta \in \mathcal{Z_{A}} \cap I_{p-j+2} \right\}
					=
		\left\{ U^{j-1}\alpha \mid 
		\alpha \in \mathcal{Z_{A}} \cap I_{j} \right\}
					\).

Since \( k \) is odd and
\( Q_{\alpha} = -Q_{\alpha^{\prime}} \) we have
\begin{equation*}
    \sum_{\alpha \in \mathcal{Z}_{\mathcal{A}}} Q_{\alpha}(z)^{-k}
    = - \sum_{\alpha \in \mathcal{Z}_{\mathcal{A}}} 
      Q_{\alpha^{\prime}}(z)^{-k}.
\end{equation*}
Using this to replace half of 
\( q^{*} \) in \eqref{RPFsymmetricodd1*}
we have
\begin{equation*}
    q^{*}(z) = \sum_{\ell=1}^{L} \frac{d_{\ell}}{2}
           \sum_{\alpha \in \mathcal{Z}_{\mathcal{A}_{\ell}}} 
       \left( Q_{\alpha}(z)^{-k} - Q_{\alpha^{\prime}}(z)^{-k} \right).
\end{equation*}
Now every 
\( \alpha \in \mathcal{Z}_{\mathcal{A}_{\ell}} \) is positive
and thus in one of the intervals 
\( I_{j} \), \( 2 \leq j \leq p \), 
of \eqref{linedecomp}.
We write \( q^{*} \) in a way that exhibits these intervals as
\begin{equation}
    \label{RPFsymmetricodd3*}
    q^{*}(z) = \sum_{\ell=1}^{L} \frac{d_{\ell}}{2}
           \sum_{j=2}^{p}
           \sum_{\alpha \in \mathcal{Z}_{\mathcal{A}_{\ell}} \bigcap I_{j}} 
	       \left( Q_{\alpha}(z)^{-k} - Q_{\alpha^{\prime}}(z)^{-k} \right).
\end{equation} 
Now  
\( 2 \leq j \leq p \) if and only if
\( 2 \leq p-j+2 \leq p \), so
\begin{equation*}
\begin{split}
	\sum_{\alpha \in \mathcal{Z}_{\mathcal{A}_{\ell}} \bigcap I_{j}} 
    Q_{\alpha^{\prime}}(z)^{-k}
    &= 
    \sum_{\alpha \in \mathcal{Z}_{\mathcal{A}_{\ell}} \bigcap I_{p-j+2}} 
    Q_{\alpha^{\prime}}(z)^{-k} \\
    &= 
    \sum_{\alpha \in \mathcal{Z}_{\mathcal{A}_{\ell}} \bigcap I_{j}} 
    Q_{U^{j-1}\alpha}(z)^{-k}.
\end{split}
\end{equation*}
We have used Lemma \ref{lemma:mappingpoles} for the second equality.
We relabel the constant for each irreducible system of poles and 
combine this with 
\eqref{RPFsymmetricodd3*} to write
\begin{equation}
    \label{RPFsymmetricodd5*}
    q^{*}(z) = \sum_{\ell=1}^{L} c_{\ell}
           \sum_{j=2}^{p}
           \sum_{\alpha \in \mathcal{Z}_{\mathcal{A}_{\ell}} \bigcap I_{j}} 
           \left( Q_{\alpha}(z)^{-k} - Q_{U^{j-1}\alpha}(z)^{-k} \right).
\end{equation}

We turn our attention to the 
second relation for the remainder term.
We define 
\begin{equation}
    \label{remaindertermpartdef}
    R(s;a,b) 
    = (a-b)^{k}\int_{0}^{\infty}\frac{y^{s}}{(iy-a)^{k}(iy-b)^{k}} \frac{dy}{y},
\end{equation}
for \( s \) a complex variable with \( 0 < \re(s) < 2k \) and
for \( a \) and \( b \) nonzero real numbers.
We substitute \eqref{RPFsymmetricodd5*}
into \eqref{remainderterm2}
and use \eqref{RPFterm}
to write an alternative expression for the remainder term,
\begin{equation}
    \label{remainderterm5}
    R(s) = -\sum_{\ell=1}^{L} c_{\ell} D_{\ell}^{-k/2}
           \sum_{j=2}^{p}
           \sum_{\alpha \in \mathcal{Z}_{\mathcal{A}_{\ell}} \bigcap I_{j}} 
           \left( R\left( s;\alpha,\alpha^{\prime} \right) 
                - R\left( s;U^{j-1}\alpha,U^{j-1}\alpha^{\prime} \right)  \right).
\end{equation}
We define a mapping 
acting on expressions of the form
\eqref{remaindertermpartdef} by
\begin{equation}
    \label{rhomapping}
    \rho\left( R\left( s;a,b \right) \right)
    = -R\left( 2k-s;a-\lambda,b-\lambda \right).
\end{equation}
We extend the action of \( \rho \) to 
linear combinations of terms of the form 
\eqref{remaindertermpartdef} by linearity,
and note that \( \rho \) is of order \( p \).
The next lemma will allow us to 
write the second relation for \( R(s) \)
using \( \rho \).
\begin{lemma}
	\label{alternativerho}
	\( -R \left( 2k-s;a-\lambda,b-\lambda \right)
	= R \left( s;U^{-1}a,U^{-1}b \right) \)
\end{lemma}
The proof of the lemma uses the integral definition 
\eqref{remaindertermpartdef}, a change of variables, 
and some simple manipulations.

\begin{theorem}
	\label{}
	Fix \( p \geq 3 \), let \( \lambda = \lambda_{p} \), and let
	\( \rho \) be the mapping defined by \eqref{rhomapping}.
	Suppose that \( R(s) \)
	is a remainder term for the Dirichlet series 
	that corresponds to an entire modular integral 
	with a Hecke-symmetric rational period function on \( G_{p} \).
	Then \( R(s) \) satisfies the (second) relation
	\begin{equation}
    \label{RTrelation2}
		R + \rho(R) + \rho^{2}(R) + \cdots + \rho^{p-1}(R) = 0.
	\end{equation}
\end{theorem}
\begin{proof}
	Suppose that \( \hat{R} \) is an expression in the domain of \( \rho \).
	Then a calculation shows that for any positive integer \( m \)
	the expression 
	\( \hat{R} - \rho^{m}\left( \hat{R} \right) \) 
	satisfies the second relation 
	\eqref{RTrelation2}.
	We use Lemma \ref{alternativerho} to rewrite the expression
	\eqref{remainderterm5} for \( R(s) \) as 
	\begin{equation}
    \label{}
    R(s) = -\sum_{\ell=1}^{L} c_{\ell} D_{\ell}^{-k/2}
           \sum_{j=2}^{p}
           \sum_{\alpha \in \mathcal{Z}_{\mathcal{A}_{\ell}} \bigcap I_{j}} 
           \left( R\left( s;\alpha,\alpha^{\prime} \right) 
                - \rho^{p-j+1}
                \left( R\left( s;\alpha,\alpha^{\prime} \right) \right) \right).
                \nonumber
	\end{equation}
	Thus \( R(s) \) is a linear combination of terms of the form
	\begin{equation*}
	    \hat{R} - \rho^{m}\left( \hat{R} \right),
	\end{equation*}
	so \( R(s) \) satisfies \eqref{RTrelation2}.
\end{proof}

\section{The converse theorem}

We now prove the following converse to 
Theorem \ref{directtheorem}.
\begin{theorem}
	\label{conversetheorem}
	Fix \( p \geq 3 \),
	let \( \lambda = \lambda_{p} \), and
	let \( k \) be an odd positive integer.
	Suppose the Dirichlet series 
	\eqref{phi}
	converges absolutely in the half-plane
	\( \re(s) > \gamma \).
	Suppose that the function 
	\( \Phi(s) \) defined by
	\eqref{PhiDS}
	satisfies
	\newline
	(a)  
	\( \Phi(s) \) has a meromorphic continuation to the 
	\( s \)-plane with, at worst, simple poles at integer points;
	\newline
	(b)  
	\( \Phi(s) \) is bounded in every lacunary vertical strip of the form
	\eqref{lacunaryVS}
	and
	\newline
	(c)  
	\( \Phi(s) \) satisfies the functional equation
	\eqref{PhiFE}
	with 
	\( R(s) \) given by 
	\eqref{remainderterm4}.
	Then \( \phi(s) \) is the Dirichlet series associated
	with an entire automorphic integral of weight 
	\( 2k \) on \( G_{p} \) with Hecke-symmetric rational period function
	given by 
	\eqref{RPFsymmetricoddexplicit} and
	\eqref{RPFpoleatzero}.
\end{theorem}

\begin{proof}
	We write 
	\( s = \sigma + it \),
	with \( \sigma, t \in \R \).
	Since \eqref{phi}
	converges absolutely in the half-plane
	\( \sigma > \gamma \),
	we have 
	\( a_{n} = \mathcal{O}\left( n^{\gamma - 1} \right) \),
	so 
	\( F(z) = \sum _{n=0}^{\infty} a_{n} e^{2 \pi inz/\lambda} \)
	converges for 
	\( z \in \uhp \).
	We follow Riemann \cite{Rie1859} and Hecke \cite{Hec36, Hec38},
	take the inverse Mellin transform of
	\( \Phi(s) \) and
	rearrange to get
	\begin{equation}
    \label{}
	    \frac{1}{2\pi i} \int_{c-i\infty}^{c+i\infty}
	    \Phi(s)y^{-s}ds
	    = \sum_{n=1}^{\infty} a_{n}e^{-2\pi ny/\lambda}
	    = F(iy),
	    \nonumber
	\end{equation}
	for any \( c > 0 \).
	The interchange of the sum and integral is valid by 
	Stirling's formula 
	\cite[p. 224]{Cop35}
	\begin{equation}
    \label{stirlingsformula}
	    \left| \Gamma(\sigma + it) \right|
	    \sim \sqrt{2\pi} |t|^{\sigma - 1/2} e^{-\pi|t|/2},
	\end{equation}
	as \( |t| \rightarrow \infty \).
	Let \( M \) be a positive integer with 
	\( M > \gamma \) and \( M > 2k \),
	and fix \( c \) with \( M < c < M+1 \).
	We move the line of integration from 
	\( \sigma = c \) to \( \sigma = 2k - c \).
	In order to do so we use the fact that
	\begin{equation*}
	    \lim_{|T| \rightarrow \infty} 
	    \int_{2k-c+iT}^{c+iT} \Phi(s) y^{-s} ds
	    = 0,
	\end{equation*}
	which follows from the fact that 
	\( \Phi(s) \) is bounded in the lacunary vertical strip
	\( S(2k-c,c;t_{0}) \),
	along with
	Stirling's formula \eqref{stirlingsformula}
	and the 
	Phragm\'en Lindel\"of Theorem \cite[p. 180]{Tit39}.
	The integrand 
	\( \Phi(s)y^{-s} \) 
	has possible poles at the integers between the lines of integration, 
	except it
	does not have a pole at 
	\( s = M \) since 
	\( \phi(s) \) converges absolutely for 
	\( \sigma = M \).
	We pick up the residues at these poles, so
	\begin{equation}
    \label{}
	    F(iy) = \frac{1}{2\pi i} \int_{2k-c-i\infty}^{2k-c+i\infty}
			    \Phi(s)y^{-s}ds
			    + \frac{1}{2\pi i} \sum_{m=2k-M}^{M-1} 
			    \res_{s=m} \left\{ \Phi(s) y^{-s} \right\}.
			    \nonumber
	\end{equation}
	We use the functional equation 
	\eqref{PhiFE} and a change of variables to calculate that
	\begin{equation}
	    \label{FdsMR1}
	    \left( F \mid T \right)(iy) - F(iy) 
	    = - \frac{1}{2\pi i} \int_{2k-c-i\infty}^{2k-c+i\infty}
			    R(s)y^{-s}ds
			    - \frac{1}{2\pi i} \sum_{m=2k-M}^{M-1} 
			    \res_{s=m} \left\{ \Phi(s) \right\} y^{-m}.
	\end{equation}

	In order to evaluate this integral 
	we invert the formula in Lemma \ref{lemma:MTformulasym}.
	Since 
	the integral in \eqref{MTformulasym} converges absolutely for
	\( 0 < \sigma < 2k \) and since
	\( \frac{1}{(iy - \delta)^{k}(iy - \epsilon)^{k}} \)
	is of bounded variation,
	we have the inverse Mellin transform
	\cite[Theorem 9a]{Wid46}
	\begin{multline}
	    \label{invMTformulasym}
		\frac{1}{2 \pi i}
		\int_{d-i\infty}^{d+i\infty}
		\{ i^{s} \delta^{s-k} B(2k-s,s-k)
		  _{2}F_{1}\left[ k, 1-k; k-s+1;\frac{\epsilon}{\epsilon-\delta} \right] \\
		  + i^{s} \epsilon^{s-k} B(s,k-s)
		_{2}F_{1}\left[ k, 1-k; s-k+1;\frac{\epsilon}{\epsilon-\delta} \right]\}
		y^{-s} ds \\
		  = \frac{(\delta - \epsilon)^{k}}{(iy - \delta)^{k}(iy - \epsilon)^{k}},
	\end{multline}
	for \( 0 < d < 2k \) and \( y > 0 \),
	and for \( \delta<0<\epsilon \), 
	\( \epsilon<0<\delta \), 
	\( 0<\epsilon<\delta \), or
	\( \delta<\epsilon<0 \).
	We move the line of integration for the integral in 
	\eqref{FdsMR1}
	from \( \sigma = 2k-c \) to \( \sigma = 1/2 \),
	and pick up the negatives of the residues of
	\( R(s)y^{-s} \) at
	\( s = 2k-M, 2k-M+1, \dots, 0 \).
	We may do so because
	\begin{equation*}
	    \lim_{|T| \rightarrow \infty} 
	    \int_{2k-c+iT}^{1/2+iT} R(s) y^{-s} ds
	    = 0,
	\end{equation*}
	which follows from 
	the fact the the Beta functions in \eqref{remainderterm4} 
	have exponential decay as \( |t| \rightarrow \infty \)
	and all of the other functions are bounded as
	\( |t| \rightarrow \infty \).
	Thus
	\begin{multline}
	    \label{FdsMR2}
	    \left( F \mid T \right)(iy) - F(iy) 
	    = - \frac{1}{2\pi i} \int_{1/2-i\infty}^{1/2+i\infty}
			    R(s)y^{-s}ds 
			    + \frac{1}{2\pi i} \sum_{m=2k-M}^{0} 
			    \res_{s=m} \left\{ R(s) \right\} y^{-m}  \\
			    - \frac{1}{2\pi i} \sum_{m=2k-M}^{M-1} 
			    \res_{s=m} \left\{ \Phi(s) \right\} y^{-m}.
	\end{multline}
	We substitute \eqref{remainderterm4} into the integral
	and
	apply \eqref{invMTformulasym} to get
	\begin{multline}
	    \label{FdsMR3}
	    \left( F \mid T \right)(iy) - F(iy) 
	    = \sum_{\ell=1}^{L} d_{\ell} D_{\ell}^{-k/2}
           \sum_{\alpha \in \mathcal{Z}_{\mathcal{A}_{\ell}}} 
	       \frac{(\alpha - \alpha^{\prime})^{k}}
	      {(iy-\alpha)^{k}(iy-\alpha^{\prime})^{k}}  \\
			    + \frac{1}{2\pi i} \sum_{m=2k-M}^{0} 
			    \res_{s=m} \left\{ R(s) \right\} y^{-m}  \\
			    - \frac{1}{2\pi i} \sum_{m=2k-M}^{M-1} 
			    \res_{s=m} \left\{ \Phi(s) \right\} y^{-m}  \\
		= q^{*}(iy) + \sum_{m=2k-M}^{0} a_{m}(iy)^{-m}
			+ \sum_{m=2k-M}^{M-1} b_{m}(iy)^{-m},
			\nonumber
	\end{multline}
	where
	\( a_{m} = \frac{i^{m}}{2\pi i} \res_{s=m} \left\{ R(s) \right\} \)
	and 
	\( b_{m} = \frac{-i^{m}}{2\pi i} \res_{s=m} \left\{ \Phi(s) \right\} \).
	The identity theorem gives us
	\begin{equation*}
	    \left( F \mid T \right)(z) - F(z) = q(z),
	\end{equation*}
	for \( z \in \uhp \), where 
	\begin{equation*}
	    q(z) = q^{*}(z) + \sum_{m=2k-M}^{0} a_{m}z^{-m}
			+ \sum_{m=2k-M}^{M-1} b_{m}z^{-m},
	\end{equation*}
	is a rational function of \( z \).
	Thus \( F \) is an entire automorphic integral of weight
	\( 2k \) on \( G_{p} \)
	with RPF \( q \).
	Now \( q^{*} \) is a Hecke-symmetric RPF, so
	\( q(z) - q^{*}(z) \) must be an RPF with a pole only at zero.
	Thus \( q(z) - q^{*}(z) \) must have the form \eqref{RPFpoleatzero}.
\end{proof}

\begin{remarks}
	(\emph{i})  The fact that \( q(z) - q^{*}(z) \) has the form given in
	\eqref{RPFpoleatzero}
	means that 
	\( b_{m} = 0 \) for \( 2k < m < M \).
	This implies that
	\( \Phi(s) \) cannot have a pole at any of the values
	\( s = 2k+1, 2k+2, \dots, M-1 \).
	This is an additional restriction 
	since we originally assumed only
	that \( \phi(s) \)
	converges absolutely in some half-plane
	\( \re(s) > \gamma \).

	(\emph{ii})  The fact that \( q(z) - q^{*}(z) \) has the form given in
	\eqref{RPFpoleatzero}
	also means that 
	\( a_{m} + b_{m} = 0 \) for 
	\( 2k-M \leq m < 2k \) with \( m \neq 0 \)
	(and \( m \neq 1 \) if \( 2k = 2 \)).
	This implies that for these values of 
	\( m \)
	\begin{equation}
	    \label{rescancel}
	    \res_{s=m} \left\{ R(s) - \Phi(s) \right\} =0.
	\end{equation}
	The fact that this holds for \( 2k-M \leq m <0 \) 
	implies that even though
	\( \Phi(s) \) and \( R(s) \) both have poles at negative integer values,
	the residues must all cancel.

	(\emph{iii})  
	If we apply the functional equation \eqref{PhiFE} 
	to \eqref{rescancel} we get
	\begin{equation*}
	    \res_{s=m} \left\{ \Phi(2k-s) \right\} =0,
	\end{equation*}
	for 
	\( 2k-M \leq m < 2k \) with \( m \neq 0 \)
	(and \( m \neq 1 \) if \( 2k = 2 \)),
	or
	\begin{equation*}
	    \res_{s=n} \left\{ \Phi(n) \right\} =0,
	\end{equation*}
	for
	\( 0 < n \leq M \) with \( n \neq 2k \)
	(and \( n \neq 1 \) if \( 2k=2 \)).
	For \( 2k < n < M \) this is our observation in 
	(\emph{i}).
\end{remarks}

\end{document}